\documentclass{amsart}

\usepackage[xurl]{fi-preamble}
\usepackage{subfiles}

\author{Kaya Arro}
\title{\FIt-Calculus and Representation Stability}
\date{June 21, 2023}

\begin{document}

\begin{abstract}
    We introduce a functor calculus for functors $\FI\to\Vc$, which we call $\FI$-objects, for $\FI$ the category of finite sets and injections and $\Vc$ a stable presentable $\infty$-category. We show that $n$-homogeneous $\FI$-objects are classified by representations of $\Sys{n}$ in $\Vc$, allowing us to associate ``Taylor coefficients'' to an $\FI$-object. We show that these Taylor coefficients, in aggregate, themselves carry the structure of an $\FI$-object, and we show that, up to the vanishing of certain Tate constructions, ``analytic'' $\FI$-objects can be recovered from their $\FI$-object of Taylor coefficients.

    We then establish a close relationship between our $\FI$-calculus and the phenomenon of representation stability for $\FI$-modules, suggesting that $\FI$-calculus be understood as the extension of representation stability to the $\infty$-categorical setting. In this context, we show how representation-theoretic information about a representation stable $\FI$-module can be read off from its $\FI$-module of Taylor coefficients.
\end{abstract}

\maketitle

\tableofcontents

\section{Introduction}

\subsection{Outline}

We describe a flavor of functor calculus for functors $\FI\to\Vc$, which we call $\FI$-objects, where $\FI$ is the category of finite sets and injections and $\Vc$ is a stable presentable $\infty$-category. At the outset, our study into this $\FI$-calculus was inspired by an analogy to Weiss' orthogonal calculus, but we soon realized that $\FI$-calculus is in fact a homotopical extension of the ideas of representation stability described in, among many other notable works, \cite{fimodstab}.

In \Cref{sectaytow} we define, for $n\in\N$, $n$-polynomial $\FI$-objects as those sending certain $n+1$-cubes to limit cubes, prove that our definition of an $n$-polynomial $\FI$-object is equivalent to the criterion that the $\FI$-object be ``presented in degree at most $n$,'' (a characterization analogous to one for representation stable $\FI$-modules) and note that every $\FI$-object admits a universal approximation by an $n$-polynomial $\FI$-object, giving rise to a Taylor tower. We introduce the $\infty$-category of ``formal Taylor towers,'' which we call ``formal Taylor towers,'' define an analytic $\FI$-object as one which is an iterated limit of polynomial $\FI$-objects, define a category of ``convergent formal Taylor towers,'' and show that the $\infty$-categories of analytic $\FI$-objects and of convergent formal Taylor towers are equivalent:

\begin{restatable*}{theorem}{anlyconv}\label{anlyconv}
	The Taylor tower functor $\Pf$ determines an equivalence
	\[
		\Pf:\FIVAnly\simeq\ExSeqVConv:\lim
	\]
\end{restatable*}

In \Cref{seccoeff}, we define an $n$-homogeneous $\FI$-object as an $n$-polynomial $\FI$-object whose universal $n-1$-polynomial approximation vanishes. We show, in analogy to other functor calculi, that the $\infty$-category of $n$-homogeneous $\FI$-objects is equivalent to the $\infty$-category of $\Sy_n$-objects. For an $\FI$-object $E$, we call the $\Sy_n$-object corresponding to the $n$-homogeneous layers of the Taylor tower of $E$ the $n$th Taylor coefficient of $E$. A priori, the Taylor coefficients of $E$ together form a symmetric sequence; we show that this symmetric sequence extends to an $\FI$-object. We define an operation $\Delf$ on $\FI$-objects which behaves like a derivative and use it to show that the Taylor tower of a finitely supported $\FI$-object is trivial, a result we summarize with the slogan, ``an analytic $\FI$-object is determined by its germ at infinity.''

We then provide conditions under which a formal Taylor tower (and hence in particular an analytic $\FI$-object) can be recovered from its $\FI$-object of Taylor coefficients; these results generalize \cite[Theorem~2.5.1]{samsnowden} which has been generalized in a different fashion as \cite[Theorem~B]{fitails}.

More specifically, we define ``tame'' formal Taylor towers -- more general than finitely presented $\FI$-objects -- and show that the aggregate Taylor coefficient functor restricts to an equivalence on tame formal Taylor towers:

\begin{restatable*}{corollary}{exseqcoeffequiv}\label{exseqcoeffequiv}
	We have an equivalence
	\[
		\Cf:\ExSeqVTame\simeq\FIVcoTame
	\]
\end{restatable*}

We formulate a yet weaker condition we call ``self-tameness'' which still ensures that a formal Taylor tower can be recovered from its Taylor coefficients.

\begin{restatable*}{corollary}{selftameequiv}\label{selftameequiv}
	We have an equivalence
	\[
		\core\Cf:\core\ExSeqVTame\simeq\core\FIV\selfcoTame
	\]
\end{restatable*}

We observe that for many choices of $\Vc$ of interest -- for example any $\Q$-linear $\infty$-category, and in particular $\QSp$, the $\infty$-category of rational chain complexes -- \emph{all} formal Taylor towers are tame, so that in such contexts a formal Taylor tower can always be recovered from its Taylor coefficients.

In \Cref{secrational} we concern ourselves with the case when $\Vc=\QSp$ and seek to show that representation stability for $\FI$-modules is an emanation of $\FI$-calculus. We have the following theorem:

\begin{restatable*}{theorem}{excrepstab}\label{excrepstab}
    For some $n\in\N$, let $E$ be an $n$-polynomial $\FI$-object taking values in rational spectra with finitely generated homology groups. Then the $\FI$-modules $H_i(E)$ are representation stable.
\end{restatable*}

Combining this result with \Cref{exseqcoeffequiv}, we calculate an explicit dictionary allowing us to read off the representations appearing in the stable part of a representation-stable $\FI$-module from its coefficient $\FI$-module by proceeding one homogeneous layer at a time:

\begin{restatable*}{corollary}{coeffrep}\label{coeffrep}
    For $E$ an $n$-homogeneous rational $\FI$-object with $\Cf_n E\cong V(\mu)$ for $\mu\vdash n$, $E|_{\FIs{\geq 2n}}\cong V(\mu)_\bullet|_{\FIs{\geq 2n}}$.
\end{restatable*}

Our results therefore suggest that a larger family of rational $\FI$-objects -- the analytic ones -- deserve consideration under the mantle of representation stability even when their homology $\FI$-modules fail to be representation stable, since their behavior is nonetheless controlled by the same functor calculus phenomena and is still recorded by their Taylor coefficient $\FI$-objects.

Additionally, because $n$-polynomial $\FI$-objects which are eventually concentrated in a particular homological dimension may not be concentrated in that dimension on sufficiently small sets, we observe that $\FI$-calculus illuminates the existence of ``good'' pre-stable behavior involving the interaction of homology groups in different dimensions in the pre-stable range. This suggests searching for such good behavior in the pre-stable ranges of real-world $\FI$-objects of interest as a way of extending downwards the lower bounds on their ``good'' behavior.



\subsection{Review of presentable stable \inftyt-categories}

An $\infty$-category $\cat{J}$ is \emph{finite} if its classifying space $B\cat{J}$ is equivalent to a finite CW-complex. A stable $\infty$-category $\Vc$ is one which is both finitely complete and finitely cocomplete and such that for $\cat{I},\cat{J}$ finite $\infty$-categories and for every functor
\begin{align*}
    F:\cat{I}\times\cat{J}\to\Vc
\end{align*}
the canonical morphism
\[
    \colim_{j\in\cat{J}}\lim_{i\in\cat{I}}F(i,j)\to\lim_{i\in\cat{I}}\colim_{j\in\cat{J}}F(i,j)
\]
is invertible. This characterization, stated in the framework of derivators, is due to Moritz Rahn and Michael Shulman and proven in \cite{stable}. We similarly have that for every
\[
    G:\cat{I}\to\Vc
\]
and for every $X\in\Vc$, the canonical morphisms
\begin{align*}
    \colim_{i\in\cat{I}}\Vc\qty(X,G(i))
        &\to\Vc\qty(X,\colim_{i\in\cat{I}} G(i))\\
    \Vc\qty(\lim_{i\in\cat{I}}G(i),X)
        &\to\lim_{i\in\cat{I}}\Vc\qty(G(i),X)
\end{align*}
are invertible.

The traditional definition, found, for example, in \cite{HA}, is that an $\infty$-category $\Vc$ is stable if it is finitely complete, finitely cocomplete, the canonical morphism from the initial object to the terminal object is invertible (i.e. $\Vc$ is pointed), and fiber and cofiber squares coincide, i.e. that a square diagram
\[\begin{tikzcd}
    A\ar[r]\ar[d]
        &B\ar[d]\\
    0\ar[r]
        &C
\end{tikzcd}\]
is a pullback if and only if it is a pushout. Note that a morphism in a stable $\infty$-category is an isomorphism if and only if either its fiber or its cofiber is contractible (or, equivalently, both its fiber and its cofiber are contractible).


Note that given two objects $X,Y\in\Vc$ in a stable $\infty$-category, the morphism object $\Vc(X,Y)$ carries the structure of a spectrum. For further details, see \cite[Chapter 1]{HA}.

The condition that $\Vc$ be presentable plays no explicit role in the paper except to ensure that $\Vc$ admit the necessary limits and colimits and to allow us to use the adjoint functor theorem for presentable $\infty$-categories. It is satisfied in all examples of interest.

We reserve the term ``sub-$\infty$-category'' for full sub-$\infty$-categories. We call a sub-$\infty$-category of an arbitrary $\infty$-category \emph{reflective} if the inclusion functor admits a left adjoint, which we call the \emph{reflection} functor. The dual notion is called a \emph{coreflective} sub-$\infty$-category. Given $\Vc'\subseteq\Vc$ a reflective sub-$\infty$-category of a pointed presentable $\infty$-category $\Vc$, the collection of objects $X\in\Vc$ such that for all $Y\in\Vc'$
    \[
        \Vc\qty(X,Y)\cong0
    \]
is called the \emph{left orthogonal complement} of $\Vc'$ and is a coreflective sub-$\infty$-category of $\Vc$. The dual notion is called the \emph{right orthogonal complement} of a reflective sub-$\infty$-category. If $\Vc$ is stable, left and right orthogonal complement are, up to equivalence, inverse operations. For $S$ a class of objects of $\Vc$, we call the intersection of all (co)reflective sub-$\infty$-categories of $\Vc$ containing the objects of $S$ the (co)reflective sub-$\infty$-category of $\Vc$ \emph{generated} by $S$.

\subsection{Notation}

Throughout, $\Spa$ refers to the $\infty$-category of spaces, $\Sp$ to that of spectra, and $\Sph$ to the sphere spectrum. Fix $\Vc$ a stable presentable $\infty$-category. Given $X\in\Vc$ and an $\infty$-groupoid $Y\in \Spa$, denote by $Y\otimes X$ the colimit of the functor $Y\to\Vc$ that is constantly $X$. The functor
\[
    -\otimes X:\Spa\to\Vc
\]
is a left adjoint so it extends canonically along $\Sigma^\infty_+$ to a left adjoint
\[
    -\otimes X:\Sp\to\Vc
\]
from $\Sp$, the stabilization of $\Spa$, and we use the same notation to describe this extended functor. Right adjoint to $-\otimes X$ is the enrichment of $\Vc$ in $\Sp$: $\Vc(X,-):\Vc\to\Sp$.

Dually, for $Y\in\Spa$ and $X\in\Vc$, we denote by $Y\power X$ the limit of the functor $Y\to\Vc$ that is constantly $X$, so that we have a functor
\[
    -\power X:\Spa\op\to \Vc
\]
We call a functor $\FI\to\Vc$ an \emph{$\FI$-object} of $\Vc$ and we denote the category of such functors $\FIV$. We sometimes conflate a set and its cardinality; e.g. when we compare two sets with symbols such as $\le$, we are really comparing their cardinalities. For $\Phi$ a property, we denote by $\FIs{\Phi}$ the full subcategory of $\FI$ spanned by those sets satisfying $\Phi$; $\FIs{\leq n}$ is a typical example. We write $\FIsV{\Phi}$ for the evident functor category.

When two $\infty$-categories are canonically equivalent up to an insignificant level of ambiguity, we sometimes conflate them. As an example, we denote by $\Sys{n}$ both the category with sole object the set \(\{1,\ldots,n\}\) and morphisms bijections and the subcategory of $\FI$ spanned by all sets with cardinality $n$.

For $\cat{C}$, $\cat{D}$ small $\infty$-categories and $\cat{E}$ a presentable $\infty$-category, when there is a canonical functor $\cat{C}\to\cat{D}$, we write
\[
    \Lans{\cat{C}}{\cat{D}}:Fun\qty(\cat{C},\cat{E})\to Fun\qty(\cat{D},\cat{E})
\]
for left Kan extension, leaving the functor $\cat{C}\to\cat{D}$ implicit. Similarly, we denote right Kan extension by
\[
    \Rans{\cat{C}}{\cat{D}}:Fun(\cat{C},\cat{E})\to Fun(\cat{D},\cat{E})
\]

\section{Polynomial \FIt-objects and Taylor towers}\label{sectaytow}

\subsection{Polynomial \FIt-objects}

\begin{definition}\label{excision}
	For $n\in\FI$, we define the \emph{$n$-cube category} to be $\FIs{/n}$, equivalently the powerset lattice of $n$. We define a \emph{standard cube} to be a diagram in $\FI$ determined by finite sets $S\subseteq S'$ in which the vertices are sets $T$ such that $S\subseteq T\subseteq S'$ and the morphisms are the inclusions. We say that a standard cube determined by $S\subseteq S'$ is a \emph{standard $n$-cube} if $S'\setminus S\cong n$.
	We say that an $\FI$-object $E$ is \emph{$n$-polynomial} if it sends each standard $n+1$-cube to a limit diagram (often called a Cartesian cube) in $\Vc$. We denote the $\infty$-category of $n$-polynomial $\FI$-objects in $\Vc$ with the notation $\ExcsV{n}$.
\end{definition}

\begin{remark}\label{standardsuffices}
	Call an $n$-cube $J:\FIs{/n}\to\FI$ \emph{semi-standard} if there exist finite sets $S,T$ and a function $f:T\to n$ such that $J(x)=S\sqcup f^{-1}(x)$ for all $x\subseteq n$. We view these semi-standard cubes as homologs of strongly coCartesian cubes of Goodwillie calculus. An $\FI$-object is $n$-polynomial if and only if it sends all semi-standard $n+1$-cubes to Cartesian cubes. We do not use this fact, so we omit the proof in the interest of brevity.
\end{remark}

\begin{recollection}\label{cubes}
	Recall that the \emph{total fiber} of an $n$-cube $J:\FIs{/n}\to\cat{C}$ for $\cat{C}$ a pointed $\infty$-category is the fiber of the canonical morphism
	\[
		J(\emptyset)\to \lim_{\emptyset\neq S\subseteq n}J(S)
	\]
	We denote the total fiber of $J$ by $\tofib J$ or $\tofib_{S\subseteq n}J(S)$. Recall that the dual notion is called the \emph{total cofiber} of $J$.
	Recall also that for $J:\FIs{/k\sqcup m}\to\cat{C}$ we have an isomorphism
	\[
		\tofib_{S\subseteq k\sqcup m}J(S)\cong \tofib_{T\subseteq k}\left(\tofib_{T'\subseteq m}J(T\sqcup T')\right)
	\]
 	Recall that when $f:X\to Y$ is a morphism in a stable $\infty$-category $\cat{C}$, $\cofib f\cong \Sigma \fib f$, so given an $n$-cube $J:\FIs{/n}\to\cat{C}$, we can regard $n$ as the disjoint union of its singleton subsets and apply the preceding result repeatedly to obtain that $\tocofib J\cong \Sigma^n\tofib J$. For more details, see \cite[Proposition~5.5.4]{cubebook}.
	Since in a stable $\infty$-category a morphism is an isomorphism if and only if its fiber is contractible, we have that in a stable $\infty$-category, an $n$-cube is Cartesian if and only if it is coCartesian. We therefore could have defined $n$-polynomial $\FI$-objects to be those sending semi-standard (or standard) $n+1$-cubes to coCartesian cubes, and we will make use of this characterization.
\end{recollection}

\begin{corollary}\label{nestedexc}
	For $m\geq n$, $\Exc_n\Vc\subseteq\Exc_{m}\Vc$.
	\begin{proof}
		For $E\in\Exc_n\Vc$ it is enough to verify that the total fiber of the image under $E$ of any standard $m+1$-cube is 0. By Recollection~\ref{cubes}, this is equivalent to the total fiber of an $(m-n)$-cube of total fibers of the images of standard $n+1$-cubes, each of which is $0$ by assumption.
	\end{proof}
\end{corollary}

\begin{proposition}\label{rightapprox}
	The full sub-$\infty$-category $\Exc_n\Vc$ of $\FIV$ is reflective. We denote its reflection functor $\Pf_n$.
\begin{proof}
	Limits commute with Cartesian cubes and filtered colimits commute with coCartesian cubes, which are Cartesian cubes because $\Vc$ is stable. It follows that $\Exc_n\Vc$ is closed under limits and filtered colimits in $\FIV$. The result follows from the adjoint functor theorem \cite[Corollary~5.5.2.9]{HTT}.\end{proof}
\end{proposition}

\begin{definitions}\label{analytic}
	We say that an $\FI$-object is \emph{polynomial} if it is $n$-polynomial for some $n\in\N$. We denote by $\FIVAnly$ the reflective subcategory of $\FIV$ generated by the polynomial $\FI$-objects and we call its objects \emph{analytic}.
\end{definitions}

\begin{definition}\label{repdef}
	Given $n\in\FI$ and $X\in\Vc$, we call $\FI$-objects isomorphic to those of the form $\FI(n,-)\otimes X$ \emph{representable}. For brevity, we denote $\Fs{n,X}\defeq \FI(n,-)\otimes X$. Recall that all $\FI$-objects are iterated colimits of representable $\FI$-objects. When $\Vc=\Sp$ and $X=\Sph$, we simply write $F_n$ for $\Fs{n,\Sph}$. 
\end{definition}

\begin{proposition}\label{polynomialrep}
	We have that for all $X\in\Vc$ and $n\in\FI$, $\Fs{n,X}\in\ExcsV{n}$.
	\begin{proof}
		Because $-\otimes X:\Spa\to\Vc$ is a left adjoint, it suffices to show that the functors
		\[
			\FI(S,-):\FI\to\Spa
		\]
		send standard $n+1$-cubes to coCartesian $n+1$-cubes.
		
		We will use the following fact. Let $f:\FIs{/n}\to\FIs{/m}$ be a functor which preserves meets (limits, intersections) and let $g:\FIs{/m}\to\Spa$ be the functor sending the subsets of $m$ to themselves understood as discrete spaces. Then $gf$ is a coCartesian cube if
		\[
			f(n)=\bigcup_{i\in n}f\left(n\setminus \{i\}\right) 
		\]

		For $S\subseteq T\subseteq S'$ and $S\subseteq T'\subseteq S'$,
		\[
			\FI(n,T\cap T')\cong \FI(n,T)\cap\FI(n,T')
		\]
		Each subset of $S'$ of cardinality $n$ is a subset of $S'\setminus\{i\}$ for some $i\in S'\setminus S$ exactly when $n<|S'\setminus S|$ -- i.e. when the standard cube in question is a standard $m$-cube for $m>n$.
	\end{proof}
\end{proposition}

\begin{theorem}\label{leftapprox}
	We have an equivalence of categories
	\[
		\Lans{\FIs{\leq n}}{\FI}:\FIsV{\leq n}\simeq\ExcsV{n}:\Ress{\FI}{\FIs{\leq n}}
	\]
	\begin{proof}
		For $n,k\in\N$, denote by $\Exc_{n,\leq k}\Vc$ the full sub-$\infty$-category of $\FIs{\leq k}\Vc$ spanned by functors sending all standard $n+1$-cubes in $\FIs{\leq k}$ to Cartesian $n+1$-cubes in $\Vc$. Because $\FI\cong\colim_{k\in\N}\FIs{\leq k}$, we have
		\[
			\FIV\cong \lim_{k\in\N\op}\FIs{\leq k}\Vc
		\]
		where the inverse limit is taken over the restriction functors $\Ress{\FIs{\leq k+1}}{\FIs{\leq k}}$, and because every standard $n+1$-cube in $\FI$ lies in $\FIs{\leq k}$ for some $k\in\N$, we also have that
		\[
			\ExcsV{n}\cong \lim_{k\in\N\op}\Exc_{n,\leq k}\Vc
		\]
		Because $\FIs{\leq k}\to\FIs{k+1}$ is fully faithful, $\Lans{\FIs{\leq k}}{\FIs{\leq k+1}}$ is right inverse to $\Ress{\FIs{\leq k+1}}{\FIs{\leq k}}$.

		For the other composition, let $n\geq k$ and $E\in \Exc_{n,\leq k+1}\Vc$. The counit
		\[
			\colim_{S\subsetneq k+1}E(S)\cong\Lans{\FIs{\leq k}}{\FIs{\leq k+1}} \Ress{\FIs{\leq k+1}}{\FIs{\leq k}}E(k+1)\overset{\varepsilon_{k+1}}{\to}E(k+1)\cong \colim_{S\subsetneq k+1}E(S)
		\]
		is an isomorphism, with the last isomorphism following from Recollection~\ref{cubes}.
	\end{proof}
\end{theorem}

\begin{observation}\label{qncolim}
	It follows that $\ExcsV{n}$ is a coreflective sub-$\infty$-category of $\FIV$ (an easier way to see the result is the adjoint functor theorem). We denote the coreflection functor $\Qf_n$. Note that
	\[
		\Qf_n E\cong \Lans{\FIs{\leq n}}{\FI}\Ress{\FI}{\FIs{\leq n}}E
	\]
	Note that for all $E\in\FIV$,
	\[
		E\cong \colim_{n\in\N}\Qf_n E
	\]
	since for all $m\geq n$, $\Qf_m E(n)\to E(n)$ is an isomorphism.
\end{observation}


\subsection{Formal Taylor towers}

\begin{definition}\label{exseq}
	We define the $\infty$-category $\ExSeqV$ of \emph{formal Taylor towers} in $\Vc$ to be
	\[
		\ExSeqV\defeq\lim\cdots\to\FIV\overset{\Pf_{1}}{\to}\FIV\overset{\Pf_{0}}{\to}\FIV
	\]
 so that a formal Taylor tower is a collection of $\FI$-objects $\{E_i\}_{i\in\N}$ (but we will often simply denote a given formal Taylor tower with a single capital Latin letter, such as $E$) equipped with, for $m\geq n$, compatible isomorphisms $\Pf_n E_m\cong E_n$. The reflection morphisms give us a tower 
	\[\begin{tikzcd}
		\vdots\ar[d]\\
		E_n\ar[d]\\
		\vdots\ar[d]\\
		E_0
	\end{tikzcd}\]
	We have a functor $\Pf:\FIV\to \ExSeqV$ given by
	\[
		(\Pf E)_i\defeq \Pf_i E
	\]
	For $E\in\FIV$, we call $\Pf E$ the \emph{Taylor tower} of $E$.
\end{definition}

\begin{definitions}\label{exseqs}
	We denote by $\ExSeqsV{n}$ the full sub-$\infty$-category of $\ExSeqV$ spanned by formal Taylor towers $\{E_i\}$ such that for all $m\geq n$, the morphism $E_m\to E_n$ is an isomorphism. We denote by
	\[
		\Jf_n:\ExSeqV\to\ExSeqsV{n}
	\]
	the coreflection functor. We say that a formal Taylor tower is \emph{convergent} if it is a colimit in $\ExSeqV$ of a diagram taking values in $\bigcup_{n\in\N}\ExSeqsV{n}$. We denote by $\ExSeqVConv$ the full sub-$\infty$-category of $\ExSeqV$ spanned by convergent formal Taylor towers.
\end{definitions}

\begin{lemma}\label{JisQ}
	For $E\in\ExSeqV$,
	\[
		\lim_{i\in\N\op}(\Jf_n E)_i\cong \Qf_n \qty(\lim_{i\in\N\op} E_i)
	\]
\end{lemma}

\begin{lemma}\label{towers}
	Let $\cat{E}$ be an arbitrary presentable $\infty$-category and let
	\[
		\cat{E}_0\subseteq\cat{E}_1\subseteq\cdots\subseteq\cat{E}_n\subseteq\cdots\subseteq\cat{E}
	\]
	be an increasing sequence of presentable reflective sub-$\infty$-categories with reflection functors $L_n$. Denote by $\cat{E}_\infty$ the full sub-$\infty$-category of $\cat{E}$ spanned by limits in $\cat{E}$ of diagrams taking values in $\bigcup_{n\in\N}\cat{E}_n$. Then $\cat{E}_\infty$ is reflective with reflection functor
	\[
		L_\infty\cong \lim_{n\in\N\op} L_n
	\]
	When $\cat{E}$ is stable, we obtain a dual theorem for presentable coreflective sub-$\infty$-categories by taking left orthogonal complements.
	\begin{proof}
		Let $\cat{E}_\infty'$ denote the intersection of all presentable reflective sub-$\infty$-categories of $\cat{E}$ which contain the union of the $\cat{E}_i$ and with reflection functor $L_\infty'$. Then $\cat{E}'_\infty$ is a presentable reflective sub-$\infty$-category of $\cat{E}$ and is the closure under iterated limits of
		\[
			\bigcup_{i\in\N}\cat{E}_i
		\]
		It follows that the canonical natural transformation $L_\infty'\qty(\id\to L_\infty)$ is an isomorphism. Then because for any $E\in\cat{E}$, $L_\infty E\in\cat{E}_\infty'$, it follows that $L_\infty\cong L_\infty'$ and hence $\cat{E}_\infty'=\cat{E}_\infty$.
	\end{proof}
\end{lemma}

\begin{corollary}\label{anlyconvcoref}
	$\FIVAnly\subseteq\FIV$ and $\ExSeqVConv\subseteq\ExSeqV$ are reflective and coreflective respectively. We denote their reflection and coreflection functors $\Pf_\infty$ and $\Jf_\infty$ respectively.
\end{corollary}

\anlyconv
\begin{proof}
	First, note that because
	\[
		\Pf:\FIV\to\ExSeqV
	\]
	factors through $\ExSeqVConv$, so does its right adjoint
	\[
		\lim_{i\in\N\op}:\ExSeqV\to\FIV
	\]
	so every analytic $\FI$-object is the limit of a convergent formal Taylor tower.
	Let
	\[
		E\in\ExSeqVConv
	\]
	We have
	\begin{align*}
		\Pf\left(\lim_{i\in\N\op} E_i\right) 
			& \cong \Pf\left(\colim_{n \in \N} \Qf_n \lim_{i \in \N\op} E_i\right)\\
			& \cong \Pf\left(\colim_{n \in \N} \lim_{i \in \N\op} \left(\Jf_n E\right)_i\right)\\
			& \cong \colim_{n \in \N} \Pf\left(\lim_{i \in N} \left(\Jf_n E\right)_i\right)\\
			& \cong \colim_{n \in \N} \left(\Jf_n E\right)_i\\
			& \cong E
	\end{align*}
	establishing that
	\[
		\Pf\circ\qty(\lim_{i\in\N\op}-)\simeq \id_{\ExSeqVConv}
	\]
	while
	\[
		\qty(\lim_{i\in\N\op}-)\circ\Pf\simeq\id_{\FIVAnly}
	\]
	follows from Lemma \ref{towers}.
\end{proof}

\section{Taylor coefficients}\label{seccoeff}

\subsection{Homogeneous \FIt-objects}

\begin{definition}\label{homg}
	We say that $E\in\ExcsV{n}$ is \emph{$n$-homogeneous} if $\Pf_{n-1}E=0$. We denote the full sub-$\infty$-category of $n$-homogeneous $\FI$-objects $\HomgsV{n}$. We define
	\[
		\Df_n\defeq \fib \qty(\Pf_n\to \Pf_{n-1}):\FIV\to\HomgsV{n}
	\]
	We say that $\Df_nE$ is the \emph{$n$th layer} of the Taylor tower of $E$. More generally, we can speak of the $n$th layer of any formal Taylor tower, and we denote this construction also by $\Df_n$.
\end{definition}

\begin{definition}\label{cohomg}
	We say that $E\in\FIV$ is \emph{$n$-cohomogeneous} if $E$ is in the image of $\Lans{\Sys{n}}{\FI}$. Equivalently, $E\in\FIV$ is $n$-cohomogeneous when $E\in\ExcsV{n}$ and $\Qf_{n-1}E\cong 0$. We denote the category of $n$-cohomogeneous $\FI$-objects $\coHomgsV{n}$. We define
    \[
        \Rf_n\defeq \cofib \qty(\Qf_{n-1}\to \Qf_n):\FIV\to\coHomgsV{n}
    \]
\end{definition}

\begin{proposition}\label{homogclassif}
	When restricted to $\HomgsV{n}$ and $\coHomgsV{n}$  respectively, the functors $\Rf_n$ and $\Df_n$ are inverses. In particular, $n$-homogeneous $\FI$-objects are classified by $\Sy_n$-objects.
	\begin{proof}
		Given $E\in\HomgsV{n}$, we consider the following commutative diagram.
		\[\begin{tikzcd}
			0\ar[d]\ar[r]&E\ar[d]\ar[r]&\Df_n\Rf_nE\ar[d]\\
			\Qf_{n-1}E\ar[d]\ar[r]&E\ar[d]\ar[r]&\Rf_nE\ar[d]\\
			\Qf_{n-1}E\ar[r]&0\ar[r]&\Pf_{n-1}\Rf_nE
		\end{tikzcd}\]
		We begin by considering the middle row. This is a fiber sequence. The bottom row is $\Pf_{n-1}$ applied to the middle row and therefore also a fiber sequence. The top row is the fiber of the natural transformation from the middle row to the bottom row and is therefore also a fiber sequence. This proves that
		\[
			\id_{\HomgsV{n}}\cong \Df_n\Rf_n
		\]
		The other direction follows from a similar argument.
	\end{proof}
\end{proposition}

\begin{definition}\label{taycoeffs}
	Given $E\in\FIV$ and $n\in\FI$, we define
	\[
		\Cf_nE\defeq \Rf_n\Df_nE(n)\cong \tocofib_{T\subseteq n}\Df_nE(T)
	\]
 	This is a $\Sy_n$-object, and we call $\Cf_n E$ the \emph{$n$th Taylor coefficient} of $E$.
\end{definition}

\subsection{The aggregate coefficient functor}

The $\Cf_n$ are left adjoint functors and so controlled by their restrictions to representable $\FI$-objects. We therefore wish to calculate $\Cf_n\Fs{S,X}$, and the first step toward that goal is calculating $\Pf_n\Fs{S,X}$.

\newcommand{\limcs}[3]{\lim_{\overset{#3 \subseteq #2}{|#3|\leq #1}}}

\begin{proposition}\label{coreppn}
	For $S\in\FI$ and $X\in\Vc$,
	\[
		\Ps{n}\Fs{S,X}\cong \limcs{n}{S}{T}\Fs{T,X}
	\]
	\begin{proof}
		For just this proof, let us denote 
		\[
			\Fs{S,X}^{(n)}\defeq \limcs{n}{S}{T}\Fs{T,X}
		\]
		Because $\Fs{X,S}^{(n)}$ is a limit of $n$-polynomial $\FI$-objects and therefore $n$-polynomial, it is enough to show that
		\[
			\Pf_n \Fs{S,X}\cong \Pf_n \Fs{S,X}^{(n)}
		\]
		By the Yoneda lemma, it is enough to show that for all $E\in\ExcsV{n}$,
		\[
			\FIV\qty(\Fs{S,X},E)\cong \FIV\qty(\Fs{S,X}^{(n)},E)
		\]
		We have
		\begin{align*}
			\FIV\qty(\Fs{S,X},E)
				&\cong\Vc\qty(X,E(S))\\
				&\cong\Vc\qty(X,\colim_{\overset{T\subseteq S}{T\leq n}}E(T))\\
				&\cong\colim_{\overset{T\subseteq S}{T\leq n}}\FIV\qty(\Fs{T,X},E)\\
				&\cong\FIV\qty(\Fs{S,X}^{(n)},E)\qedhere
		\end{align*}
	\end{proof}
\end{proposition}

\begin{definition}
	In what follows, we denote $\Gs{n,X}\defeq \Df_n \Fs{n,X}$ and $\Gs{n}\defeq\Df_n\Fs{n}$.
\end{definition}

\begin{corollary}\label{gn}
	\[
		\Gs{n,X}\cong \tofib_{S\subseteq n} \Fs{S,X}
	\]
\end{corollary}

\begin{corollary}
	\[
		\Df_n\Fs{m,X}\cong \prod_{\overset{S\subseteq m}{|S|=n}} \Gs{S,X}
	\]
	\begin{proof}
		Using \Cref{coreppn},
		\begin{align*}
			\Ds{n}\Fs{m,X}
				&\cong\fib\Ps{n}\Fs{m,X}\to\Ps{n-1}\Fs{m,X}\\
				&\cong\fib\limcs{n}{m}{S}\Fs{S,X}\to\Ps{n-1}\limcs{n}{m}{S}\Fs{S,X}\\
				&\cong\limcs{n}{m}{S}\fib\Fs{S,X}\to\Ps{n-1}\Fs{S,X}\\
				&\cong\Rans{\qty{S\subset m:|S|=n}}{\qty{S\subset m:|S|\leq n}}\Gs{S,X}\\
				&\cong\prod_{\overset{S\subseteq m}{|S|=n}} \Gs{S,X}\qedhere
		\end{align*}
	\end{proof}
\end{corollary}

\begin{proposition}\label{repcoeffs}
	The Taylor coefficients of representable $\FI$-objects are given by
	\[
		\Cf_n\Fs{m,X}\cong \FI(n,m)\pitchfork X
	\]
	\begin{proof} 
		Using the preceding corollaries, we have
		\begin{align}
			\Cf_n \Fs{m,X}&=\tocofib_{T\subseteq n}\prod_{\overset{U\subseteq m}{|U|=n}} \tofib_{S\subseteq U}\Fs{S,X}(T)\nonumber\\
				&\cong\prod_{\overset{U\subseteq m}{|U|=n}}\tofib_{S\subseteq U}\tocofib_{T\subseteq n}\Fs{S,X}(T)\nonumber\\
				&\cong\prod_{\overset{U\subseteq m}{|U|=n}}\tocofib_{T\subseteq n}\Fs{U,X}(T)\label{eqnnotofib}\\
				&\cong\prod_{\overset{U\subseteq m}{|U|=n}}\Fs{U,X}(n)\label{eqnnotocofib}\\
				&\cong\prod_{\overset{U\subseteq m}{|U|=n}}\FI(U,n)\otimes X\nonumber\\
				&\cong\prod_{\overset{U\subseteq m}{|U|=n}}\FI(n,U)\pitchfork X\label{colimtolim}\\
				&\cong\FI(n,m)\pitchfork X\nonumber
		\end{align}
		where \cref{eqnnotofib} uses that when $|S|<n$, $\Fs{S,X}$ is $n-1$-polynomial so that
		\[
			\tocofib_{T\subseteq n}\Fs{S,X}(T)\cong 0
		\]
		and hence
		\[
			\tofib_{S\subseteq U}\tocofib_{T\subseteq n}\Fs{S,X}(T)\cong \fib\left(\tocofib_{T\subseteq n}\Fs{U,X}(T)\to 0\right)\cong \tocofib_{T\subseteq n}\Fs{U,X}(T)
		\]
		and \cref{eqnnotocofib} uses that when $|T|<|U|$,
		\[
			\Fs{U,X}(T)=\FI(U,T)\otimes X=\emptyset \otimes X\cong 0
		\]
		so that
		\[
			\tocofib_{T\subseteq n}\Fs{U,X}(T)\cong\cofib\left(0\to \Fs{U,X}(n)\right)\cong \Fs{U,X}(n)\qedhere
		\]
	\end{proof}
\end{proposition}

\begin{corollary}\label{aggregatecoeffs}
	For $E\in \FIV$, $\Cf_n E$ is functorial in $n\in\FI$, so that we obtain an aggregate Taylor coefficient functor
	\[
		\Cf:\FIV\to\FIV
	\]
	and we extend the construction to $E\in\ExSeqV$ with the formula
	\[
		\Cf:E\mapsto\lim_{n\in\N\op}\Cf\Ps{n}E
	\]
\end{corollary}

\begin{definition}
	Denote by $\SuppsV{n}$ the image of $\Rans{\FIs{\leq n}}{\FI}$. Denote the reflection functor by
	\[
		\Us{n}:\FIV\to\SuppsV{n}
	\]
\end{definition}

\begin{observation}\label{pnun}
	For $E\in\ExSeqV$,
	\[
		\Cf\Ps{n}E\cong\Us{n}\Cf E
	\]
\end{observation}

\subsection{Derivatives}

Let us give a more direct description of $\Cf E$ in terms of the $\FI$-object $E$.

\begin{notation}
	For $E\in\FIV$ and $n\in\FI$, define a new $\FI$-object
	\[
		\Delf^n E\defeq\tocofib_{S\subseteq n}E(S\sqcup -)
	\]

	Given a map $f:n\to n'$ in $\FI$, abbreviate
	\[
		n'\setminus f\defeq n'\setminus \img f
	\]
	Further, given $k\in \FI$, define
	\[
		j^\dagger::k\to n'\sqcup k
	\]
	by
	\[
		j^\dagger(a)\defeq
		\begin{cases}
			j^{-1}(a)&a\in j(n'\setminus f)\\
			a&a\notin j(n'\setminus f)
		\end{cases}
	\]
	and define
	\[
		g_{f,k}\defeq\sum_{j:n'\setminus f\to k}E\left(f\sqcup j^\dagger\right):E(n\sqcup k)\to E(n'\sqcup k)
	\]
\end{notation}

\begin{theorem}\label{whatisc}
	The Taylor coefficients of $E$ are given by
	\[
		\Cf E(n)\cong\colim_{k\in\FI}\Delf^nE(k)
	\]
	The morphism $\Cf E(f)$ is determined by the maps $g_{f,k}$.
	\begin{proof}
		Since our constructions preserve colimits in $E$, we need only verify that the theorem holds for representable $\FI$-objects, and this in turn allows us to reduce further to the case $\Vc=\Sp$ and thence to just the $F_n$.
		To verify the formula for objects, observe that
		\[
			\colim_{k\in\FI}F_m(n+k)
		\]
		can be identified with the suspension spectrum of the set of partial bijections from $m$ to $n$. The construction $\Delf^n$ commutes with suspension, and by the fact we used in the proof of \Cref{polynomialrep}, $\Delf^n$ kills off exactly those partially defined injections which do not cover $n$. A partially defined injection $m\to n$ which covers $n$ is the same data as an injection $n\to m$, and since this set is finite, its suspension spectrum is isomorphic to its dual.
		For the remainder of the proof, we make use of the isomorphism $\Cf F_m(n)\cong \Sigma^\infty \FI(n,m)_+$. We that
		\[
			\Cf F_m(f):\Cf F_m(n)\to \Cf F_m(n')
		\]
		is determined by adjointness by the map
		\[
			\FI(n,m)\mapsto\Omega^\infty\Sigma^\infty\FI(n',m)_+
		\]
		given by
		\[
			\qty(i:n\to m)\mapsto \sum_{\overset{i':n'\to m}{i'=fi}}\eta\qty(i')
		\]
		where $\eta$ is the unit of the $\Sigma^\infty_+\dashv\Omega^\infty$ adjunction and the sum is taken with respect to the $E_\infty$-structure of the infinite loop space.
		To verify the formula for morphisms, we must verify the commutativity of the following square:
		\[\begin{tikzcd}
			F_m(n+k)\ar[d,"g_{f,k}"]\ar[r]&\Cf F_m(n)\ar[d,"\Cf F_m(f)"]\\
			F_m(n'+k)\ar[r]&\Cf F_m(n')
		\end{tikzcd}\]
		By adjointness, the canonical map $F_m(n+k)\to \Cf F_m(n)$ is determined by the map 
		\[
			\FI(m,n+k)\to\Omega^\infty\Cf F_m(n)
		\]
		which sends an injection $i:m\to n+k$ to $0$ if $n\not\subseteq i(m)$ and to $\eta(i^*)$ otherwise, where
		\[
			i^*\defeq \qty(a\mapsto i^{-1}(a)):n\to m
		\]
		We can conclude by observing that given an injection $i:m\to n+k$ representing an injection $i^*:n\to m$, each injection $n'\to m$ which restricts to $i^*$ is represented once by a summand of $g_{f,k}(i)$.%
	\end{proof}
\end{theorem}

\begin{corollary}\label{delfcoeff}
	By \Cref{whatisc} and Recollection~\ref{cubes},
	\[
		\Cf \qty(\Delf^nE)(-)\cong \Cf (E)(-\sqcup n)
	\]
\end{corollary}


\begin{definition}\label{tors}
	Call an $\FI$-object $E$ \emph{finitely supported} if there exists $n\in\N$ such that $E\in\Supp_n\Vc$. Denote by $\FIVTors$ -- the $\infty$-category of \emph{torsion} $\FI$-objects -- the coreflective sub-$\infty$-category of $\FIV$ generated by the finitely supported $\FI$-objects.
\end{definition}

\begin{corollary}\label{coeffaretors}
	For any $E\in\FIV$, $\Cf E\in\FIVTors$.
	\begin{proof}
		This follows from the facts that $\Cf E$ is finitely supported when $E$ is representable, that $\Cf$ is a left adjoint functor, and that every $\FI$-object is an iterated colimit of representable $\FI$-objects.
	\end{proof}
\end{corollary}

We formalize the notion that an analytic $\FI$-object is determined by its ``germ at infinity.''

\begin{corollary}\label{torsvanish}
	For $E \in \FIVTors$, $\Pf E \cong 0$.
	\begin{proof}
		It is sufficient to prove that $\Cf E\cong0$ for $E$ finitely supported. We begin by showing that $\Cf E(0)\cong \colim E\cong 0$ when $E|_{\FIs{\geq n}}\cong 0$. By \cite[Proposition~4.1.3.1]{HTT} -- Quillen's Theorem A for quasicategories, originally due to Joyal -- it suffices to show that for all $m\in\FI$,
        \[
            B(m\downarrow \FIs{\geq n})\cong *
        \]
        where $B:\cat{C}at_\infty\to\Spa$ is the classifying space functor, since this implies that
        \[
            0=\colim_{\FIs{\geq n}}0
            \cong\colim_{\FIs{\geq n}}E|_{\FIs{\geq n}}
            \cong \colim_\FI E
        \]
        Note that we have an equivalence
        \[
            (k,f:m\to k)\mapsto k\setminus f:m\downarrow \FIs{\geq n}\cong \FIs{\geq n-m}
        \]
        so it will be sufficient to establish that $B\FIs{\geq n}\cong *$ for all $n\in\N$.
		Let
        \[
            \iota_{\geq n}:\FIs{\geq n}\to\FI\] denote the inclusion functor and write \[\kappa_{\geq n}\defeq S\mapsto S\sqcup n:\FI\to \FIs{\geq n}
        \]
        Then we have natural transformations
        \[
            \id_{\FI}\to\iota_{\geq n}\kappa_{\geq n}
        \]
        and
        \[
            \id_{\FIs{\geq n}}\to\kappa_{\geq n}\iota_{\geq n}
        \]
        each given by the canonical inclusion $S\to S\sqcup n$. Upon taking classifying spaces, these natural transformations become homotopies, so that $B\iota_{\geq n}$ and $B\kappa_{\geq n}$ are inverse homotopy equivalences between $B\FI$ and $B\FIs{\geq n}$. But $B\FI\cong *$ because $\FI$ has an initial object.
	    For any $S\in\FI$, we have that if $E$ vanishes on $\FIs{\geq n}$, then so does $\Delf^S E$, so for $E$ finitely supported,
        \[
            \Cf E(S)\cong\Cf(\Delf^S E)(0)\cong 0\qedhere
        \]
	\end{proof}
\end{corollary}

\begin{corollary}
    For $E\in\FIVAnly$ and any $n\in\N$,
    \[
        E\cong\Rans{\FIs{\geq n}}{\FI}\Ress{\FI}{\FIs{n}}E
    \]
\end{corollary}

\begin{conjecture}\label{torsanlyorthog}
    Corollary \ref{torsvanish} implies that $\FIVAnly$ is a full sub-$\infty$-category of the right orthogonal complement of $\FIVTors$. We conjecture the converse: that $\FIVAnly$ in fact is the right orthogonal complement of $\FIVTors$.
\end{conjecture}

\subsection{Recovering formal Taylor towers from coefficients}

Without developing a functor calculus, Sam and Snowden prove in \cite[Theorem~2.5.1]{samsnowden} that tails of representation stable $\FI$-modules are fully described by torsion $\FI$-modules. In this section, we integrate that result into our $\FI$-calculus and generalize it to settings other than the rational one.

\begin{notation}
	Denote by $\Zf$ the right adjoint to $\Cf$. We have used $\Cf$ to refer to functors with various domains including $\FIV$, $\ExSeqV$, and $\ExcsV{n}$, and we shall similarly use $\Zf$ to refer to functors with these various codomains.
\end{notation}

\begin{observation}
	For $E\in\FIV$,
	\begin{align*}
		\Zf E
			&\cong\int_{m\in\FI}\Cf F_m\power E(m)\\
			&\cong\int_{m\in\FI}\Sp\qty(\FI(m,-),\Sph)\power E(m)\\
			&\cong\int_{m\in\FI}\FI(m,-)\otimes E(m)
	\end{align*}
\end{observation}

\begin{lemma}\label{unitrepiso}
	The unit
	\[
		\eta_{\Fs{X,n}}:\Fs{X,n}\to\Zf\Cf\Fs{X,n}
	\]
	is an isomorphism.
	\begin{proof}
		\begin{align*}
			\Zf\Cf\Fs{X,n}
				&\cong\int_{m\in\FI}\FI(m,-)\otimes\FI(m,n)\power X\\
				&\cong\int_{m\in\FI}\FI(m,n)\power\FI(m,-)\otimes X\\
				&\cong\FI(n,-)\otimes X\\
				&\cong\Fs{X,n}
		\end{align*}
	\end{proof}
\end{lemma}

\begin{definition}\label{tame}
	We call a $\Sys{n}$-object $A$ \emph{tame} if the norm map 
	\[
		\FI(-,n)\power_{\Sys{n}} A\to\FI(-,n)\power^{\Sys{n}} A
	\]
	is a natural isomorphism. In this case we also call the cohomogeneous $\FI$-object $\Lans{\Sys{n}}{\FI}A$ tame. We also call any $m$-polynomial $\FI$-object that is a finite colimit of tame cohomogeneous $\FI$-objects tame, any formal Taylor tower of tame polynomial $\FI$-objects tame, and any analytic $\FI$-object with a tame Taylor tower tame. We denote the $\infty$-categories of such $\SysVTame{n}$, $\coHomgsVTame{n}$, $\ExcsVTame{m}$, $\ExSeqVTame$, and $\FIVAnlyTame$ respectively.

	We call an $\FI$-object \emph{cotame} if it lies in the image of $\Ress{\ExSeqV}{\ExSeqVTame}\Cf$. We denote the $\infty$-category of cotame $\FI$-objects $\FIVcoTame$. We denote by $\SuppsVcoTame{n}$ the full sub-$\infty$-category of $\SuppsV{n}$ spanned by cotame objects.
\end{definition}

\begin{theorem}\label{cohomgunit}
	Let $A\in\SysVTame{n}$ and denote
	\[
		E\defeq \FI(n,-)\otimes A
	\]
	and endow $E$ with the diagonal $\Sys{n}$-action so that
	\[
		E_{\Sys{n}}\cong\Lans{\Sys{n}}{\FI}\in\coHomgsVTame{n}
	\]
	Then
	\[
		\eta:E_{\Sys{n}}\to\Zf\Cf E_{\Sys{n}}
	\]
	is an isomorphism.
	\begin{proof}
		Fixing total orders on $n$ and $k$ determines an isomorphism of $\Sys{n}$-spaces
		\[
			\FI(n,k)\cong\Lans{*}{\Sys{n}}\qty{S\subseteq k:|S|=n}
		\]
		exhibiting $\FI(n,k)$ as a free $\Sys{n}$-space. Tensoring an $\Sys{n}$-object with a free $\Sys{n}$-space yields a free $\Sys{n}$-object, so $E(k)$ is a free $\Sys{n}$-object for each $k\in\FI$. Then by \cite[Example~6.1.6.26]{HA}, the norm maps
		\[
			\Nm{E(k)}:E(k)_{\Sys{n}}\to E(k)^{\Sys{n}}
		\]
		are isomorphisms and hence the norm map
		\[
			\Nm{E}:E_{\Sys{n}}\to E^{\Sys{n}}
		\]
		is an isomorphism.

		Because left adjoint functors of stable $\infty$-categories are exact, we have a commutative square
		\[\begin{tikzcd}[column sep=huge]
			\Cf\qty(E_{\Sys{n}})\ar[r,"\Cf\qty(\Nm{E})"]\ar[d]
				&\Cf\qty(E^{\Sys{n}})\ar[d,"f"]\\
			\qty(\Cf E)_{\Sys{n}}\ar[r,"\Nm{\Cf E}"]
				&\qty(\Cf E)^{\Sys{n}}
		\end{tikzcd}\]
		where the left arrow is an isomorphism because $\Cf$ is a left adjoint, the top arrow is $\Cf$ applied to the norm and hence an isomorphism by the preceding argument, and the bottom arrow is an isomorphism by our assumption on $A$, so the right arrow $f$ must also be an isomorphism.

		The adjunction $\Cf\dashv\Zf$ now gives us a new commutative square
		\[\begin{tikzcd}
			E^{\Sys{n}}\ar[r,"\qty(\eta_E)^{\Sys{n}}"]\ar[d,"\eta_{E^{\Sys{n}}}"]
				&\qty(\Zf\Cf E)^{\Sys{n}}\ar[d]\\
			\Zf\Cf\qty(E^{\Sys{n}})\ar[r,"\Zf f"]
				&\Zf\qty(\Cf E)^{\Sys{n}}
		\end{tikzcd}\]
		We know that all the morphisms in this square except the left one are isomorphisms, so that one must be as well.

		We play our game one last time. The naturality of the unit gives us the commutative square
		\[\begin{tikzcd}[column sep=huge]
			E_{\Sys{n}}\ar[r,"\Nm{E}"]\ar[d,"\eta_{E_{\Sys{n}}}"]
				& E^{\Sys{n}}\ar[d]\\
			\Zf\Cf\qty(E_{\Sys{n}})\ar[r,"\Zf\Cf\qty(\Nm{E})"]
				&\Zf\Cf\qty(E^{\Sys{n}})
		\end{tikzcd}\]
		and we conclude that $\eta_{E_{\Sys{n}}}$ must be an isomorphism since all the other morphisms in the square are isomorphisms.
	\end{proof}
\end{theorem}

\begin{remark}
	The bulk of the foregoing proof can be encapsulated by the claim that the following square commutes:
	\[\begin{tikzcd}
		E_{\Sys{n}}\ar[r]\ar[d]
			& E^{\Sys{n}}\ar[d]\\
		\Zf\Cf\qty(E_{\Sys{n}})\ar[r]
			&\qty(\Zf\Cf E)^{\Sys{n}}
	\end{tikzcd}\]
\end{remark}

\begin{lemma}\label{stablefinitedense}
	Suppose that $\cat{C}$ and $\cat{D}$ are stable $\infty$-catories, that $\cat{C}_0\subseteq\cat{C}$ and $\cat{D}_0\subseteq\cat{D}$ are full sub-$\infty$-categories, that each object of $\cat{C}$ and $\cat{D}$ can be expressed as a colimit of objects in $\cat{C}_0$ and $\cat{D}_0$ respectively, and that there exists an exact functor
	\[
		L:\cat{C}\to\cat{D}
	\]
	which restricts to an equivalence
	\[
		\Ress{\cat{C}}{\cat{C}_0}L:\cat{C}_0\simeq\cat{D}_0
	\]
	Then $L$ is an equivalence
	\[
		L:\cat{C}\simeq\cat{D}
	\]
	\begin{proof}
		We show that $L$ is surjective and fully faithful. Denote by
		\[
			R:\cat{D}_0\to\cat{C}_0		
		\]
		the inverse of $\Ress{\cat{C}}{\cat{C}_0}L:\cat{C}_0$. For establish surjectivity, we have that for $a\in\cat{D}$ there exists some finite $\infty$-category $\cat{I}$ and diagram
		\[
			A:\cat{I}\to\cat{D}_0
		\]
		such that
		\begin{align*}
			a
				&\cong\colim_{i\in\cat{I}}A(i)\\
				&\cong\colim_{i\in\cat{I}}LRA(i)\\
				&\cong L\qty(\colim_{i\in\cat{I}}RA(i))
		\end{align*}
		Next, for $b,c\in\cat{C}$, there exist finite $\infty$-categories $\cat{J}$ and $\cat{K}$ and functors
		\begin{align*}
			B:\cat{J}&\to\cat{C_0}\\
			C:\cat{K}&\to\cat{C_0}
		\end{align*}
		such that $b\cong\colim B$ and $c\cong\colim C$ so that we have
		\begin{align*}
			\cat{D}\qty(Lb,Lc)
				&\cong \cat{D}\qty(L\qty(\colim{j\in\cat{J}} B(j)),L\qty(\colim_{k\in\cat{K}} C(k)))\\
				&\cong \cat{D}\qty(\colim_{j\in\cat{J}} LB(j),\colim_{k\in\cat{K}} LC(k))\\
				&\cong \lim_{j\in\cat{J}}\colim_{k\in\cat{K}}\cat{D}_0\qty(LB,LC)\\
				&\cong \lim_{j\in\cat{J}}\colim_{k\in\cat{K}}\cat{C}_0\qty(B,C)\\
				&\cong \cat{C}\qty(\colim_{j\in\cat{J}} B(j),\colim_{k\in\cat{K}} C(k))\\
				&\cong \cat{C}\qty(b,c)
		\end{align*}
		establishing full faithfulness.
	\end{proof}
\end{lemma}

\begin{corollary}
	We have an equivalence
	\[
		\Cf:\ExcsVTame{n}\simeq\SuppsVcoTame{n}
	\]
\end{corollary}

\exseqcoeffequiv

\begin{corollary}\label{anlyequiv}
	The equivalence of \Cref{exseqcoeffequiv} restricts to an equivalence
	\[
		\Cf:\FIVAnlyTame\simeq\FIVTorscoTame
	\]
	\begin{proof}
		This follows from \Cref{coeffaretors} and \Cref{anlyconv}.
	\end{proof}
\end{corollary}

\begin{example}
	When $\Vc$ is $\Q$-linear, all $\Sys{n}$-objects are tame, so we have equivalences
	\begin{align*}
		\Cf:\ExSeqV&\simeq\FIV\\
		\Cf:\FIVAnly&\simeq\FIVTors
	\end{align*}
\end{example}

\begin{example}
	When $\Vc$ is the $\infty$-category of $K(m)$-local spectra for some $m\in\N$ and prime $p$, all $\Sys{n}$-objects are tame, so we have equivalences
	\begin{align*}
		\Cf:\ExSeqV&\simeq\FIV\\
		\Cf:\FIVAnly&\simeq\FIVTors
	\end{align*}
\end{example}

\begin{definition}
	For $E\in\ExSeqV$, we say that $E$ is \emph{self-tame} if for all pairs $n,m\in\N$, the map
	\[
		\FIV\qty(\Qs{m}E,\Qs{n}E)\to\FIV\qty(\Qs{m}E,\Zf\Cf\Qs{n}E)
	\]
	induced by the unit (i.e. the norm map) is an isomorphism. We denote the full sub-$\infty$-category of such $\ExSeqV\selfTame$. We say that an $\FI$-object is \emph{self-cotame} if it lies in the image of
	\[
		\Cf:\ExSeqV\selfTame\to\FIV
	\]
	and we denote the full sub-$\infty$-category of such $\FIV\selfcoTame$.
\end{definition}

\begin{notation}
	Recall that the \emph{core} of an $\infty$-category $\cat{C}$, denoted $\core\cat{C}$, is its coreflection into $\Spa$; in other words, $\core\cat{C}$ is the maximal sub-$\infty$-groupoid of $\cat{C}$.
\end{notation}

\selftameequiv

\section{Representation stability}\label{secrational}

\subsection{A bouquet}

\begin{notation}\label{uSp}
    We define
    \[
        \uSp:\Spa\to\Spa_*
    \]
    by
    \[
        \uSp: X\mapsto \cofib\qty(\qty(X\to *)_+)
    \]
    In words, we make from an unbased space $X$ a based space $\uSp X$ by taking the unreduced suspension $SX$ of $X$ and designating one of the cone points the basepoint.
\end{notation}

\begin{theorem}\label{gncon}
    For $k\ge 2n-1$, $G_n(k)$ is a wedge of copies of $\Sph$.
    \begin{proof}
        Recall that by \Cref{coreppn}, $G_n$ is the total fiber of the $n$-cube given by $F_S(k)$ as $S$ ranges over the subsets of $n$. This is equivalent to the $n$-fold desuspension of the total cofiber of the same $n$-cube, so we could equivalently show that this total cofiber is a wedge of copies of $\Sph^n$. For this it would suffice to show that the total cofiber, which we denote by $L(n,k)$, of the $n$-cube $\FI(S,k)_+$ is a wedge of copies of $S^n$, where we form the colimit in the $\infty$-category of pointed spaces rather than of spectra.

        As a first step, we will show that
        \begin{equation}\label{cofibseq1}
            \cofib\qty(L(n,k)\to L(n,k\sqcup 1))\cong\bigvee_{x\in n}\Sigma L(n\setminus \{x\},k)
        \end{equation}
        Note that
        \begin{equation}\label{cofibseq2}
            \cofib\qty(\FI(S,k)_+\to\FI(S,k\sqcup 1)_+)\cong \bigvee_{x\in S}\FI(S\setminus\{x\},k)_+
        \end{equation}
        For $x\in n$, define the $n$-cube $A_x$ by
        \[
            A_x:S\mapsto\begin{cases}\FI(S\setminus\{x\},k)_+&x\in S\\{*}&x\notin S\end{cases}
        \]
        and note that
        \[
            \tocofib_{S\subseteq n} A_x(S)\cong \cofib \qty(L\qty(n\setminus\{x\},k)\to *)\cong\Sigma L\qty(n\setminus\{x\},k)
        \]
        We can rewrite isomorphism~\ref{cofibseq2} as
        \[
            \cofib\qty(\FI(S,k)_+\to\FI(S,k\sqcup 1)_+)\cong \bigvee_{x\in n} A_x(S)
        \]
        and taking total cofibers over $S\subseteq n$ yields isomorphism~\ref{cofibseq1}.

        Next, observe that because $\FI(\emptyset,k)\cong *$, we have that
        \begin{equation}\label{Kissus}
            L(n,k)\cong \uSp \colim_{\emptyset\neq S\subseteq n}\FI(S,k)
        \end{equation}
        Let us consider the domain of the left fibration classified by the functor $\FI(-,k):\FIs{/n,>0}\to\cat{S}$. This is the partially ordered set $P(n,k)$ of tuples $(S,T,\phi)$ where $\emptyset\neq S\subseteq n$, $\emptyset\neq T\subseteq k$, and $\phi:S\cong T$ with order given by $(S,T,\phi)\leq (S',T',\phi')$ if $S\subseteq S'$, $T\subseteq T'$, and $\phi'|_S=\phi$. There is an evident isomorphism $P(n,k)\cong P(k,n)$, and because $P(n,k)$ is the domain of the left fibration classified by $\FI(-,k):\FIs{/n,>0}\to\cat{S}$, we have $NP(n,k)\cong \colim_{\emptyset\neq S\subseteq n}\FI(S,k)$ where $N$ denotes the nerve of the poset.

        The symmetry $P(n,k)\cong P(k,n)$ reveals the symmetry $L(n,k)\cong L(k,n)$ by \cref{Kissus}. Combining this with \cref{cofibseq1}, we have
        \begin{equation}\label{cofibseq3}
            \cofib\qty(L(n,k)\to L(n\sqcup 1,k))\cong\bigvee_{x\in k}\Sigma L(n,k\setminus \{x\})
        \end{equation}
        We make the inductive hypothesis that for some $n$ there exists $C(n)$ such that for all $k\geq C(n)$, $L(n,k)$ is a wedge of copies of $S^n$. Let $k\geq C(n)$ and consider the long exact sequence in homology induced by \cref{cofibseq3} (but replace $k$ with $k\sqcup 1$). By our inductive hypothesis, for $i\neq n$,
        \[
            H_{i+1}\qty(\bigvee_{x\in k\sqcup 1}\Sigma L(n,\qty(k\sqcup 1)\setminus\{x\}))\cong 0\cong H_i(L(n,k\sqcup 1))
        \]
        so necessarily
        \(
            H_i(L(n\sqcup 1,k\sqcup 1))\cong 0
        \)
        whenever $i\notin\{n,n+1\}$ and the morphism
        \begin{equation}\label{surj1}\begin{tikzcd}
                H_n\qty(L(n,k\sqcup 1))\ar[r,two heads]&H_n\qty(L(n\sqcup 1,k\sqcup 1))
        \end{tikzcd}\end{equation}
        is a surjection. Similarly, using \cref{cofibseq1} (and replacing $n$ and $k$ with $n\sqcup 1$ and $k\sqcup 1$ respectively), we have that
        \begin{equation}\label{surj2}\begin{tikzcd}
            H_n(L(n\sqcup 1,k\sqcup 1))\ar[r,two heads]&H_n (L(n\sqcup 1,k\sqcup 2))
        \end{tikzcd}\end{equation}
        is surjective. Composing morphisms~\ref{surj1} and \ref{surj2}, we have a surjection
        \begin{equation}\label{surj3}\begin{tikzcd}
            H_n(L(n,k\sqcup 1))\ar[r,two heads]&H_n(L(n\sqcup 1,k\sqcup 2))
        \end{tikzcd}\end{equation}
        We now show that the map
        \[
            \begin{tikzcd}L(n,k\sqcup 1)\ar[r]&L(n\sqcup 1,k\sqcup 2)\end{tikzcd}
        \]
        is nullhomotopic, so that morphism~\ref{surj3} and therefore also $H_n(L(n\sqcup 1,k\sqcup 2))$ must be trivial. This is easiest to see in terms of the posets $P(n,k\sqcup 1)$ and $P(n\sqcup 1,k\sqcup 2)$. Let us establish the notation $1=\{a\}$ and $2=\{a,b\}$. Consider the following (non-commutative) diagram:
        \begin{equation*}\label{nulhom}
            \begin{tikzcd}[column sep=huge,row sep=huge]
                P(n,k\sqcup 1)\ar[r,"\id"]\ar[d]\ar[rd,"f"]&P(n,k\sqcup 1)\ar[d,"g"]\\*\ar[r,"\{a\}\cong\{b\}"]&P(n\sqcup 1,k\sqcup 2)
            \end{tikzcd}
        \end{equation*}
        We let $g$ be the inclusion and
        \[
            f:(S,T,\phi)\mapsto (S\sqcup\{a\},T\sqcup\{b\},\phi\sqcup (\{a\}\cong\{b\}))
        \]
        Then there are natural transformations \(g\Rightarrow f\) and \((\{a\}\cong\{b\})\circ *\Rightarrow f\). After taking the nerve, these natural transformations become homotopies, and composing these homotopies yields a null-homotopy of $g$. But $g$ is the morphism which induces morphism~\ref{surj3}.

        We have proven the following: if there exists $C(n)$ such that for all $k\geq C(n)$, $L(n,k)$ is a wedge of $n$-spheres, then for all $k\geq C(n)$, the homology of $L(n\sqcup 1,k\sqcup 2)$ is concentrated in degree $n+1$ and is free (since by the long exact sequence from \cref{cofibseq3} it must be a subgroup of the homology of a wedge of spheres). Since $P(\emptyset,k)=\emptyset$ for all $k$, $L(\emptyset,k)=S^0$ for all $k$. For all $k\geq 1$, $L(1,k)$ is the suspension of a non-empty discrete space and therefore a wedge of circles. We will show that for $n\geq 2$ and $k\geq 2n-1$, $L(n,k)$ is the suspension of a connected space and therefore simply-connected. By induction, $L(n,k)$ is a Moore space of type $M(G,n)$ for $G$ free and is therefore a wedge of $n$-spheres.

        Let us show that when $n\geq 2$ and $k\geq 3$, $NP(n,k)$ is connected. Note that each vertex is connected by a $1$-simplex to a vertex of the form $(\{x\},\{y\},\phi)$. We will denote vertices $(\{x_1,\ldots,x_j\},\{y_1,\ldots,y_j\},\phi)$ by
        \[
            \begin{pmatrix}x_1&\phi(x_1)\\\vdots&\vdots\\x_j&\phi(x_j)\end{pmatrix}
        \]
        Let $x,x'\in n$ be distinct and $y,y',y''\in k$ be distinct. We must show that $\begin{pmatrix}x&y\end{pmatrix}$ is connected to $\begin{pmatrix}x'&y'\end{pmatrix}$, to $\begin{pmatrix}x'&y\end{pmatrix}$, and to $\begin{pmatrix}x&y'\end{pmatrix}$. We have
        \begin{equation}\label{ineq1}
            \begin{pmatrix}
                x&y
            \end{pmatrix}
            <
            \begin{pmatrix}
                x&y\\
                x'&y'
            \end{pmatrix}
            >
            \begin{pmatrix}
                x'&y'
            \end{pmatrix}
        \end{equation}
        \begin{equation}\label{ineq2}
            \begin{pmatrix}
                x&y
            \end{pmatrix}
            <
            \begin{pmatrix}
                x&y\\
                x'&y''
            \end{pmatrix}
            >
            \begin{pmatrix}
                x'&y''
            \end{pmatrix}
            <
            \begin{pmatrix}
                x&y'\\
                x'&y''
            \end{pmatrix}
            >
            \begin{pmatrix}
                x&y'
            \end{pmatrix}
        \end{equation}
        \begin{equation}\label{ineq3}
            \begin{pmatrix}
                x&y'
            \end{pmatrix}
            <
            \begin{pmatrix}
                x&y'\\
                x'&y
            \end{pmatrix}
            >
            \begin{pmatrix}
                x'&y
            \end{pmatrix}
        \end{equation}
        where we compose sequences~\ref{ineq2} and~\ref{ineq3} to obtain a path from $\begin{pmatrix}x&y\end{pmatrix}$ to $\begin{pmatrix}x'&y\end{pmatrix}$.\qedhere
    \end{proof}
\end{theorem}

\subsection{A dictionary}

Throughout the remainder of this section, we will at times treat rational vector spaces as rational spectra concentrated in dimension 0. For $X$ a spectrum, we denote by $X^\Q$ its rationalization. We use the term \emph{$\FI$-module} to describe functors
\[
    \FI\to\FinQVect
\]
where $\FinQVect$ is the category of finite-dimensional rational vector spaces.

\excrepstab
\begin{proof}
    When $n=0$, $E$ must be constant. An objectwise-finite $\FI$-module is representation stable if and only if it is finitely generated by \cite[Theorem~1.13]{fimodstab}. By \cite[Theorem~1.3]{fimodstab}, the category of $\FI$-modules is Noetherian.
    
    Note that the homology of an $m$-cohomogeneous $\FI$-object is generated entirely in degree $m$ and consider the cofiber sequence
    \[
            \Qf_{n-1}E\to E\to \Rf_n E
    \]
    Since $H_i (E)$ is the extension of a sub-$\FI$-module of $H_i \Rf_n E$ by a quotient $\FI$-module of $H_i\qty(\Qf_{n-1}E)$, both of which are finitely generated, $H_i (E)$ is itself finitely generated and hence representation stable.
\end{proof}

We know from \Cref{gncon} that the homology of $G_n(k)$ is concentrated in dimension 0 for $k\geq 2n-1$. Let us get to know $H_0(G_n,\Q)$ better. First we recall some facts about the representation theory of the symmetric groups.

\begin{recollection}\label{symreprec}
    Given an irreducible rational representation $V$ of $\Sys{n}$, its complexification $V\otimes_\Q\C$ is an irreducible complex representation of $\Sys{n}$, so the representation theory of $\Sys{n}$ is the same over any characteristic 0 field, regardless of algebraic closure. This is a well-known fact, but can be seen for instance from the facts that the so-called \emph{Specht modules} can be defined over the integers \cite[Section~2.3]{sagansymrep} and account for all irreducible complex representations \cite[Section~2.4]{sagansymrep}.

    Recall that a \emph{partition} is a finite, weakly decreasing sequence of natural numbers $\lambda=(\lambda_1,\ldots,\lambda_j)$. We write
    \[
        |\lambda|\defeq\sum_{i}\lambda_i
    \]
    and say that $\lambda$ is a \emph{partition of $|\lambda|$}. The notation $\lambda\vdash n$ is synonymous with $|\lambda|=n$.

    Recall e.g. from \cite[Section~2.3]{sagansymrep} that the Specht modules of $\Sys{n}$ are in bijection with partitions of $n$. We denote by $V(\lambda)$ the Specht module corresponding to a partition $\lambda$. Given a partition $\lambda$ and $k\geq \lambda_1+|\lambda|$, we define $\lambda[k]$ to be the partition $(k-|\lambda|,\lambda)$. We define \(w(\lambda)\defeq|\lambda|-\lambda_1\) and we say that $w(\lambda)$ is the \emph{weight of $\lambda$}. Observe that $w(\lambda[k])=|\lambda|$ when $\lambda[k]$ is defined.

    For $\mu\vdash n$, we have a $\Sys{n}$-representation $M^\mu$ called a \emph{Young permutation representation} and defined in \cite[Section~2.1]{sagansymrep}. By \cite[Section~2.10]{sagansymrep}, for $\lambda\vdash n$,
    \[
        \dim\qty(\Sys{n}\FinQVect(V(\lambda),M^\mu))=K_{\lambda,\mu}
    \]
    where $K_{\lambda,\mu}$ is a \emph{Kostka number}: the number of \emph{semistandard Young tableaux of shape $\lambda$ and content $\mu$}. This means the following. For $\lambda=(\lambda_1,\ldots,\lambda_j)\vdash n$, we consider $n$ boxes arranged in $j$ rows with $\lambda_i$ boxes in row $i$. A semistandard tableau of shape $\lambda$ and content $\mu$ is a way of filling these boxes with natural numbers such that the number $i$ occurs $\mu_i$ times, the columns of our tableau are strictly increasing from top to bottom, and the rows of our tableau are weakly increasing from left to right. A \emph{standard $\lambda$-tableau} means a semistandard tableau of shape $\lambda$ and content $(1^n)$.

    Finally, recall some notation from \cite{fimodstab}. Given a rational $\Sys{n}$-representation $V:\Sys{n}\to \QVect$, we define $M(V)_\bullet\defeq\Lans{\Sys{n}}{\FI}V$, the left Kan extension of $V$ to $\FI$. We shall also make use of the indecomposable $\FI$-module $V(\lambda)_\bullet$ defined in \cite[Proposition~3.1.4]{fimodstab}, which is representation stable and satisfies $V(\lambda)_k\cong V(\lambda[k])$ when $k\geq \lambda_1$ and $V(\lambda)_k\cong 0$ otherwise.
\end{recollection}

Observe that $\QG_n$ has a natural action of $\Sys{n}$ because $F_n$ does. In the following theorem, we characterize $\QG_n$ along with its $\Sys{n}$-action.

\begin{theorem}\label{gnweightn}
    We have an isomorphism in the category $Fun(\Sys{n}\times \FIs{\geq 2n},\Sp)$
    \[
        \QG_n\cong \bigoplus_{\lambda\vdash n}V(\lambda)\boxtimes V(\lambda)_\bullet
    \]
    \begin{proof}
        Given a rational $\Sy_k$-spectrum $X$ and $\lambda\vdash k$, define $\chi_\lambda(X)$ to be the Euler characteristic of the spectrum $\Sp_{\Sy_k}(V(\lambda),X)$ if it exists. Then it follows from \Cref{coreppn} that
        \begin{equation}\label{euchar1}
            \chi_\lambda\qty(\QG_n(k))=\sum_{0\leq i\leq n}(-1)^{n-i}\binom{n}{i}\chi_\lambda\qty(\QF_i(k))
        \end{equation}
        We observe that $\QF_i(k)\cong M^{\qty(k-i,1^i)}$ so that
        \[
            \chi_{\lambda}\qty(\QF_i(k))=K_{\lambda,\qty(k-i,1^i)}
        \]
        In a semistandard tableau of shape $\lambda$ and content $\qty(k-i,1^i)$, the $k-i$ $1$s must be in the left-most boxes of the first row of the tableau, so we need only consider $\lambda\vdash k$ such that $\lambda_1\geq k-n$, and since $k\geq 2n$, we have $k-\lambda_1\leq n\leq k-i$. This means that there are no boxes directly below any box after the $(k-i)$th box in the first row of our tableau, so for any subset of cardinality $\lambda_1-k+i$ of $\{2,\ldots,i+1\}$, there is exactly one way to fill out the rest of the first row of our tableau. Define $\lambda'\defeq (\lambda_2,\ldots,\lambda_j)$. Once we have chosen how to fill the top row of our semistandard $\lambda$-tableau under construction, any standard $\lambda'$ tableau determines a unique semistandard $\lambda$ tableau with the given first row (since all of the relevant numbers in the first row are 1s and all of our remaining numbers after the first row are unique and greater than 1). Note that $\lambda_1=k-w(\lambda)$, so we have shown that
        \begin{equation}
            \label{Kostanum}K_{\lambda,(k-i,1^i)}=\binom{i}{\lambda_1-k+i}K_{\lambda',1^{w\qty(\lambda)}}=\binom{i}{i-w(\lambda)}K_{\lambda',1^{w\qty(\lambda)}}=\binom{i}{w(\lambda)}K_{\lambda',1^{w\qty(\lambda)}}
        \end{equation}
        Combining this with \cref{euchar1}, we have
        \begin{equation}\label{euchar2}
            \chi_\lambda\qty(\QG_n(k))=K_{\lambda',1^{w\qty(\lambda)}}\sum_{0\leq i\leq n}(-1)^{n-i}\binom{n}{i}\binom{i}{w(\lambda)}
        \end{equation}
        Note that the quantity \(\binom{n}{i}\binom{i}{w(\lambda)}\) is the number of pairs $A,B$ with $n\supseteq A\supseteq B$ such that $|A|=i$ and $|B|=w(\lambda)$. We could also count these pairs by first choosing $B\subseteq n$ and then choosing $A\setminus B\subseteq n\setminus B$. This observation gives us the identity
        \[
            \binom{n}{i}\binom{i}{w(\lambda)}=\binom{n}{w(\lambda)}\binom{n-w(\lambda)}{i-w(\lambda)}
        \]
        Substituting this into \cref{euchar2}, we have
        \begin{align}
            \chi_\lambda\qty(\QG_n(k))
                &=K_{\lambda',1^{w\qty(\lambda)}}\binom{n}{w(\lambda)}\sum_{0\leq i\leq n}(-1)^{n-i}\binom{n-w(\lambda)}{i-w(\lambda)}\nonumber\\
                &=K_{\lambda',1^{w\qty(\lambda)}}\binom{n}{w(\lambda)}\sum_{0\leq j\leq n-w(\lambda)}(-1)^{n-w(\lambda)-j}\binom{n-w(\lambda)}{j}\nonumber\\
                \label{binomuse}&=K_{\lambda',1^{w\qty(\lambda)}}\binom{n}{w(\lambda)} (1-1)^{n-w(\lambda)}\\
                &=\begin{cases}
                    K_{\lambda',1^{w\qty(\lambda)}}\binom{n}{w(\lambda)}&n=w(\lambda)\\
                    0&n\neq w(\lambda)
                \end{cases}\nonumber\\
                \label{weightpart}&=\begin{cases}
                    \chi_\lambda\qty(\QF_n(k))&n=w(\lambda)\\
                    0&n\neq w(\lambda)
                \end{cases}
        \end{align}
        where \cref{binomuse} uses the binomial theorem.
        This establishes that for $\lambda\vdash k$, $\chi_\lambda\qty(\QG_n(k))\neq 0$ if and only if $w(\lambda)=n$. This implies that for $k\geq 2n$ and $V$ an irreducible $\Sy_k$-representation, \(\Sp_{\Sy_k}\qty(V,\QG_n(k))\neq 0\) if and only if $V\cong V(\lambda')_k$ for some $\lambda'\vdash n$.

        Note that by \Cref{gncon}, $H_0\qty(\QG_n)$ is a sub-$\Sys{n}\times \FIs{\geq 2n}$-module of $H_0\qty(\QF_n)$ and by \cref{weightpart} $H_0\qty(\QG_n)$ consists of exactly the weight-$n$ submodules of $H_0\qty(\QF_n)$. Observe that $\QF_n=M(\Q[\Sys{n}])$. By Maschke's Theorem,
        \begin{align*}
            \Q[\Sys{n}]
                &\cong \bigoplus_{\lambda\vdash n}\End(V(\lambda))\\
                &\cong \bigoplus_{\lambda\vdash n}V(\lambda)^*\boxtimes V(\lambda)\\
                &\cong \bigoplus_{\lambda\vdash n}V(\lambda)\boxtimes V(\lambda)
        \end{align*}
        where the last isomorphism holds because in characteristic zero, finite-dimensional representations of $\Sys{n}$ are self-dual. We therefore have that
        \[
            \QF_n\cong \bigoplus_{\lambda\vdash n}V(\lambda)\boxtimes M\qty(V(\lambda))
        \]
        By \cite[Lemma~3.2.3 and~Proposition 3.4.1]{fimodstab}, the weight-$n$ irreducible sub-representations of the $M(V(\lambda))(k)$ form a sub-$\FI$-module of $M(V(\lambda))$ and indeed are exactly $V(\lambda)_\bullet$.\qedhere
    \end{proof}
\end{theorem}

\coeffrep
\begin{proof}
    Restricting to $\FIs{\geq 2n}$, we have
    \begin{align*}
        E
            &\cong V(\mu)\smashp_{\Sys{n}}\QG_n\\
            &\cong \bigoplus_{\lambda\vdash n}\qty(V(\mu)\otimes V(\lambda))_{\Sys{n}}\boxtimes V(\lambda)_\bullet\\
            &\cong \bigoplus_{\lambda\vdash n}\qty(V(\mu)\otimes V(\lambda))^{\Sys{n}}\boxtimes V(\lambda)_\bullet\\
            &\cong \bigoplus_{\lambda\vdash n}\qty(V(\mu)^*\otimes V(\lambda))^{\Sys{n}}\boxtimes V(\lambda)_\bullet\\
            &\cong\bigoplus_{\lambda\vdash n}\mathrm{Hom}_{\Sys{n}}(V(\mu),V(\lambda))\boxtimes V(\lambda)_\bullet\\&\cong V(\mu)_\bullet\qedhere
    \end{align*}
\end{proof}

Since every rational $\Sys{n}$-spectrum is a direct sum of ((de)suspensions of) spectra of the form appearing in the hypothesis of \Cref{coeffrep}, we now have an elementary dictionary allowing us to translate between rational $\Sys{n}$-spectra and $n$-homogeneous rational $\FI$-objects -- in other words, we have made the equivalence in \Cref{homogclassif} explicit (in the rational case). In fact, we have the following additional corollary.

\begin{corollary}\label{obwisesum}
    For $E$ a rational $\FI$-object and $k\geq 2n$, there is an isomorphism of $\Sys{k}$-spectra
    \[
        \Pf_n E(k)\cong \bigvee_{i\leq n}D_i E(k)
    \]
    \begin{proof}
        The result follows from an inductive argument. The case $n=0$ holds since $\Pf_0=\Df_0$. For the inductive step, we apply \Cref{coeffrep} and Schur's Lemma to the long exact sequence in homology associated to the fiber sequence $\Df_n E(k)\to \Pf_n E(k)\to \Pf_{n-1}E(k)$ to see that the two sides of our equation agree on the level of homology (including the $\Sys{k}$-action). Our result follows from the fact that $\Q[\Sys{k}]$ is semisimple.
    \end{proof}
\end{corollary}



\section{Further questions}\label{secfurther}



Our work leaves several questions open and also suggests multiple avenues of further research.

We recall Conjecture~\ref{torsanlyorthog} that $\FIVAnly$ is the right orthogonal complement of $\FIVTors$. This would establish that every $\FI$-object is canonically an extension of its analytic part by its torsion part.

We would like to extend \Cref{selftameequiv} to give a convenient description of the fibers $\core\Cf$ in the absence of any tameness assumptions.

We are interested in investigating interactions between $\FI$-calculus and Goodwillie calculus. For instance, for $E\in\FIV$ and $E':\cat{V}\to\cat{V}'$, can $\Pf (E'\circ E)$ be recovered from $\Pf E$ and the Goodwillie tower of $E'$? Is there a chain rule relating $\Cf E$ and the Goodwillie derivatives of $E'$ to $\Cf (E'\circ E)$? Going in the opposite direction, if $E\in\FIVAnly$ and $E'\in\FI\cat{V}'^\Anly$, what can we say about the Goodwillie tower of $\Lan_{E}E'$ based on $\Pf E$ and $\Pf E'$?

When $\Vc$ is endowed with a symmetric monoidal structure in which compact objects are dualizable, we can construct symmetric monoidal structures on $\FIVTors$ that allow us to extend
\[
    \Cf:\FIV\to\FIVTors
\]
to a symmetric monoidal functor for $\FIV$ endowed with either the objectwise symmetric monoidal structure or with Day convolution (using the disjoint union symmetric monoidal structure on $\FI$). We have not examined this construction in depth, but it may prove useful as an aid in calculations in future work.

Finally, there are a wide range of categories other than $\FI$ on  which we can set up similar functor calculi in which we recover our results, appropriately modified, at least through \Cref{homogclassif} and \Cref{taycoeffs}.

Note that $\FI$ is equivalent to $\Embd{0}$, the category of $0$-manifolds and embeddings. Our techniques allow us to define functor calculi to the $\infty$-categories $\Embd{d}$. We note that these functor calculi differ from embedding calculus, which deals with contravariant rather than covariant functors on $\Embd{d}$ and in which polynomial approximations can be obtained by taking right Kan extensions, a technique which does not apply in our covariant setting.

Further, suppose we an $\infty$-category $\cat{C}$ equipped with a Cartesian fibration
\[
    \varpi:\cat{C}\to\Embd{d}
\]
Then we can extend our functor calculus from $\Embd{d}$ to $\cat{C}$ and again we retain our classification of homogeneous functors.

Returning to the case $d=0$, if we assume that $\varpi$ is also symmetric monoidal for $\FI$ equipped with the disjoint union symmetric monoidal structure, we obtain a notion of $c$-polynomial functors for any $c\in\cat{C}$.

One example of an $\infty$-category $\cat{C}$ equipped with a Cartesian fibration to $\FI$ is the braid category of a manifold $M$ -- that is, the category in which objects are finite ordered tuples of distinct points in $M$ and morphisms are isotopies, i.e. braids, from one tuple to another, with the possibility of additional points in the codomain of a morphism. Some examples of Cartesian fibrations over $\Embd{d}$ are given by $\infty$-categories of manifolds equipped with ``local'' structure, e.g. vector bundles, possibly with some structure on the fibers. It may be possible in such cases to study functors from such $\infty$-categories using the style of functor calculus developed here as well as a sort of parametrized orthogonal calculus simultaneously.



\clearpage
\printbibliography

\end{document}